\documentclass{elsart}
\usepackage{amsfonts}
\usepackage{graphicx}
\usepackage{amsmath}
\usepackage{amssymb}
\usepackage{bm}

\newtheorem{theorem}{Theorem}[section]
\newtheorem{corollary}[theorem]{Corollary}
\newtheorem{lemma}[theorem]{Lemma}

\begin{document}

\begin{frontmatter}

\title{On a relationship between the T-congruence Sylvester equation and the Lyapunov equation}

\author{Masaya Oozawa\corauthref{cor1}},
\ead{m-oozawa@na.cse.nagoya-u.ac.jp}
\author{Tomohiro Sogabe},
\author{Yuto Miyatake},
\author{Shao-Liang Zhang}
\address{Department of Computational Science and Engineering, Graduate School of engineering, Nagoya University, Furo-cho, Chikusa-ku, Nagoya 464-8603, Japan}
\corauth[cor1]{Corresponding author.}


\begin{abstract}

We consider the T-congruence Sylvester equation $AX+X^{\rm T}B=C$, where $A\in \mathbb R^{m\times n}$, $B\in \mathbb R^{n\times m}$ and $C\in \mathbb R^{m\times m}$ are given, and matrix $X \in \mathbb R^{n\times m}$ is to be determined.
The T-congruence Sylvester equation has recently attracted attention because of a relationship with palindromic eigenvalue problems. For example, necessary and sufficient conditions for the existence and uniqueness of solutions, and numerical solvers have been intensively studied. In this note, we will show that, under a certain condition, the T-congruence Sylvester equation can be transformed into the Lyapunov equation. This may lead to further properties and efficient numerical solvers by utilizing a great deal of studies on the Lyapunov equation.

\end{abstract}
\begin{keyword}
T-congruence Sylvester equation, Lyapunov equation, the tensor product
\end{keyword}
\end{frontmatter}


\section{Introduction}
\quad We consider the T-congruence Sylvester equation of the form
\begin{equation}
\label{T-con}
AX+X^{\rm T}B=C,
\end{equation}
where $A\in \mathbb R^{m\times n}$, $B\in \mathbb R^{n\times m}$ and $C\in \mathbb R^{m\times m}$ are given, and $X\in \mathbb R^{n\times m}$ is to be determined.
If $X=X^{\rm T}$, then the equation \eqref{T-con} reduces to the Sylvester equation that is widely known in control theory and numerical linear algebra.\\
\quad The T-congruence Sylvester equation \eqref{T-con} has recently attracted attention because of a relationship with palindromic eigenvalue problems \cite{unique2,palindromic}: necessary and sufficient conditions for the existence of a unique solution for every right-hand side \cite{unique1,unique2}, an algorithm to compute a unique solution of \eqref{T-con} in the case of $n=m$ \cite{general}, and a method to find the general solutions of $AX+X^\star B=O$ \cite{AX+X^TB=0}, where
$X^\star$ denotes either the transpose or the conjugate transpose of $X$. \\
\quad The standard technique to solve or to analyze the T-congruence Sylvester equation \eqref{T-con} is to find nonsingular matrices $T_1$, $T_2$, $T_3$ and $T_4$ of the form
\begin{equation}\label{T}
T_1AT_2(T_2^{-1}XT_3)+(T_1X^{\rm T}T_4^{-1})(T_4BT_3)=T_1CT_3
\end{equation}
such that the equation \eqref{T-con} is simplified. Furthermore, if the above equation \eqref{T} can be transformed into a well-known matrix equation, then one can use the studies of the corresponding matrix equation in order to solve or to analyze the original equation \eqref{T-con}. However, there seems to be no such transformation due to the existence of $X^{\rm T}$.\\
\quad In this paper, we show that, if $A$ or $B$ is nonsingular, the T-congruence Sylvester equation \eqref{T-con} can be transformed into the Lyapunov equation that is also known in control theory \cite{L_unique} or the Sylvester equation by using the tensor product, the vec operator and a permutation matrix.
This result will provide an approach to finding mathematical properties and efficient numerical solvers for the T-congruence Sylvester equation \eqref{T-con}.\\
\quad The paper is organized as follows. Section 2 describes the definition of the tensor product and the exchange of the tensor product by using permutation matrix.
Section 3 presents our main results that the T-congruence Sylvester equation can be transformed into the Lyapunov equation or the Sylvester equation.
Finally, Section 4 gives concluding remarks and future work.
\section{Preliminaries}
\quad In this section, the tensor product (also referred to as the Kronecker product) and its properties are briefly reviewed.
Let $A=[a_{ij}]\in \mathbb R^{m\times n}$ and $B\in \mathbb R^{p\times q}$, then the tensor product is defined by 
\begin{equation}\nonumber
 A\otimes B:=
 \begin{bmatrix}
  a_{11}B &a_{12}B &\cdots &a_{1n}B\\
  a_{21}B &a_{22}B &\cdots &a_{2n}B\\
  \vdots  &\vdots  &       &\vdots\\
  a_{m1}B &a_{m2}B &\cdots &a_{mn}B
 \end{bmatrix}
 \in \mathbb R^{mp\times nq}.
\end{equation}
In addition, let $C\in \mathbb R^{n\times l}$ and $D\in \mathbb R^{q\times r}$, then it follows that
\begin{equation}\nonumber
(A\otimes B)(C\otimes D)=(AC)\otimes(BD).
\end{equation}
\quad For $A=[\bm{a}_1,\bm{a}_2,\dots,\bm{a}_n]\in \mathbb R^{m\times n}$, the vec operator, vec:$\mathbb R^{m\times n}\to \mathbb R^{mn},$ is defined by 
\begin{equation}\nonumber
 {\rm vec} (A):=
 \begin{bmatrix}
  \bm{a}_1\\
  \bm{a}_2\\
  \vdots \\
  \bm{a}_n
 \end{bmatrix},
\end{equation}
and vec$^{-1}$:$\mathbb R^{mn}\to \mathbb R^{m\times n}$ is the inverse vec operater such that 
\begin{equation}\nonumber
 {\rm vec}^{-1} ({\rm vec} (A))=A.
\end{equation}
We shall use the following two lemmas to prove the main results in the next section.
\begin{lemma}$\!\!(${\rm See, e.g., \cite[p.275]{demmel})}
Let $A\in \mathbb R^{m\times n}, B\in \mathbb R^{n\times p}$. Then it follows that
\begin{align}\nonumber
{\rm vec} (AB)=(I_p\otimes A){\rm vec}(B)
              =(B^{{\rm T}}\otimes I_m){\rm vec}(A),
\end{align}
where $I_n$ denotes the $n \times n$ identity matrix.
\end{lemma}
\begin{lemma}$\!\!(${\rm \cite{tensor})}\label{per}
Let $\bm{e}_{in}$ be an $n$-dimensional column vector that has 1 in the $i$th position and 0's elsewhere, i.e.,
\begin{equation}\nonumber
 \bm{e}_{in}:=[0,0,\dots,0,1,0,\dots,0]^{{\rm T}}\in \mathbb R^n.
\end{equation}
Then for the permutation matrix
\begin{equation}\nonumber
 P_{mn}:=
  \begin{bmatrix}
  I_m\otimes \bm{e}_{1n}^{\scalebox{0.5}{\rm T}}\\
  I_m\otimes \bm{e}_{2n}^{\scalebox{0.5}{\rm T}}\\
  \vdots\\
  I_m\otimes \bm{e}_{nn}^{\scalebox{0.5}{\rm T}}
 \end{bmatrix}
 \in \mathbb R^{mn\times mn},
\end{equation}
the following properties hold {\rm:}
\begin{align}
&P_{mn}^{{\rm T}}=P_{nm},\label{first}\\
&P_{mn}^{{\rm T}}P_{mn}=P_{mn}P_{mn}^{{\rm T}}=I_{mn},\label{second}\\
&{\rm vec}(A)=P_{mn}{\rm vec}(A^{{\rm T}}), A\in \mathbb R^{m\times n},\label{third}\\
&P_{mr}(A\otimes I_r)P_{nr}^{{\rm T}}=I_r\otimes A.\label{forth}
\end{align}
\end{lemma}
\section{Main results}
In this section, we consider the T-congruence Sylvester equation \eqref{T-con} for the case $A$, $B$, $C\in \mathbb R^{n\times n}$.
\begin{theorem}\label{the1}
Let $A$, $B$, $C$, $X\in \mathbb R^{n\times n}$. Then we have, if $A$ is nonsingular, the T-congruence Sylvester equation \eqref{T-con} can be transformed into the Lyapunov equation of the form 
\begin{equation}\label{theorem_A}
\tilde{X}-M\tilde{X}M^{\rm T}=Q, 
\end{equation}
where $\tilde{X}:=AX$, $M:=B^{{\rm T}}A^{-1}$, $Q:=C-{\rm vec}^{-1}(P_{nn}{\rm vec}(MC)).$
\end{theorem}
{\bf Proof}
\quad Applying the vec operator to \eqref{T-con} and using Lemma 2.1 yield
\begin{equation}
\label{vecT}
 (I_n\otimes A){\rm vec}(X)+(B^{\rm T}\otimes I_n){\rm vec}(X^{\rm T})={\rm vec}(C).
\end{equation}
From \eqref{first}, \eqref{second} and \eqref{forth} of Lemma \ref{per}, it follows that
\begin{equation}\label{perm_exchange}
P_{nn}(B^{\rm T}\otimes I_n)P_{nn}=I_n\otimes B^{\rm T}
\Leftrightarrow (B^{\rm T}\otimes I_n)P_{nn}=P_{nn}(I_n\otimes B^{\rm T}).
\end{equation}
By using \eqref{perm_exchange} and \eqref{third} of Lemma \ref{per}, the second term of the left-hand side in \eqref{vecT} is culculated to be
\begin{align}
(B^{\rm T}\otimes I_n){\rm vec}(X^{\rm T})=(B^{\rm T}\otimes I_n)P_{nn}{\rm vec}(X)\nonumber
=P_{nn}(I_n\otimes B^{\rm T}){\rm vec}(X).\nonumber
\end{align}
Thus \eqref{T-con} is rewritten as
\begin{equation}
 \{(I_n\otimes A)+P_{nn}(I_n\otimes B^{\rm T})\}\bm{x}=\bm{c},\nonumber
\end{equation}
where $\bm{x}:={\rm vec}(X)$ and $\bm{c}:={\rm vec}(C)$.
Since $A$ is assumed to be nonsingular, it follows that $I_{n^2}=(I_n\otimes A^{-1})(I_n\otimes A)$, and thus
\begin{align}
&\{(I_n\otimes A)+P_{nn}(I_n\otimes B^{\rm T})\}(I_n\otimes A^{-1})(I_n\otimes A)\bm{x}=\bm{c}\nonumber\\
&\quad \Leftrightarrow \{I_n\otimes I_n+P_{nn}(I_n\otimes M)\}\bm{\tilde{x}}=\bm{c},\label{vecT2}
\end{align}
where $\bm{\tilde{x}}:=(I_n\otimes A)\bm{x}$.
Multiplying on the left of \eqref{vecT2} by $\{I_n\otimes I_n-P_{nn}(I_n\otimes M)\}$ and using \eqref{forth} of Lemma \ref{per} yield
\begin{align}
&\{I_n\otimes I_n-P_{nn}(I_n\otimes M)\}\{I_n\otimes I_n+P_{nn}(I_n\otimes M)\}\bm{\tilde{x}}=\bm{c}'\nonumber\\
&\quad \Leftrightarrow\{I_n\otimes I_n-P_{nn}(I_n\otimes M)P_{nn}(I_n\otimes M)\}\bm{\tilde{x}}=\bm{c}'\nonumber\\
&\quad \Leftrightarrow\{I_n\otimes I_n-(M\otimes I_n)(I_n\otimes M)\}\bm{\tilde{x}}=\bm{c}'\nonumber\\
&\quad \Leftrightarrow\{I_n\otimes I_n-(M\otimes M)\}\bm{\tilde{x}}=\bm{c}',\label{theorem1}
\end{align}
where $\bm{c}':=\{I_n\otimes I_n-P_{nn}(I_n\otimes M)\}\bm{c}$.
Applying the inverse vec operator to \eqref{theorem1}, we obtain
\begin{equation}
{\rm vec}^{-1}\{\{I_n\otimes I_n-(M\otimes M)\}\bm{\tilde{x}}\}={\rm vec}^{-1}(\bm{c}')\Leftrightarrow \tilde{X}-M\tilde{X}M^{\rm T}=Q.\nonumber
\end{equation}
Thus when $A$ is nonsingular, the T-congruence Sylvester equation \eqref{T-con} can be transformed into the Lyapunov equation.\hfill \ $\Box$
\begin{corollary}
Let $A$, $B$, $C$, $X\in \mathbb R^{n\times n}$. Then we have, if $B$ is nonsingular, the T-congruence Sylvester equation \eqref{T-con} can be transformed into the Lyapunov equation of the form
\begin{equation}\label{theorem_B}
\hat{X}-\hat{M}\hat{X}\hat{M}^{\rm T}=\hat{Q},
\end{equation}
 where $\hat{X}:=X^{\rm T}B$, $\hat{M}:=A(B^{{\rm T}})^{-1}$, $\hat{Q}:=C-{\rm vec}^{-1}(P_{nn}{\rm vec}(C\hat{M}^{\rm T}))$.
\end{corollary}
{\bf Proof}
Transposing the equation (\ref{T-con}) yields
\begin{equation}\label{the_2}
B^{\rm T}X+X^{\rm T}A^{\rm T}=C^{\rm T}.
\end{equation}
By replacing $B^{\rm T}$, $A^{\rm T}$, and $C^{\rm T}$ with $A$, $B$, and $C$ respectively, the equation \eqref{the_2} becomes \eqref{T-con}. As the result of Theorem 3.1, when B is nonsingular, the T-congruence Sylvester equation \eqref{T-con} can be transformed into the Lyapunov equation. \hfill \ $\Box$
\par
\quad As an application, if the matrix $A$ or $B$ is nonsingular, one may obtain the solution of \eqref{T-con} as follows:\par
1. solve the Lyapunov equation \eqref{theorem_A} (or \eqref{theorem_B});\par
2. solve $AX=\tilde{X}$ (or $B^{\rm T}X=\hat{X}^{\rm T}$).\par
\quad A slightly stronger condition gives a close relationship between the T-congruence Sylvester equation and the Sylvester equation as shown below.
\begin{corollary}
Let $A$, $B$, $C$, $X\in \mathbb R^{n\times n}$. Then we have, if $A$ and $B$ are nonsingular, the T-congruence Sylvester equation can be transformed into the Sylvester equation{\rm:} 
\begin{equation}
-M\tilde{X}+\tilde{X}(M^{-1})^{\rm T}=Q',
\end{equation}
where $\tilde{X}:=AX$, $M:=B^{{\rm T}}A^{-1}$, $Q':={\rm vec}^{-1}\{(M^{-1}\otimes I_n)\bm{c'}\}$.
\end{corollary}
{\bf Proof}\quad 
Multiplying on the left of \eqref{theorem1} by $(M^{-1}\otimes I_n)$ gives
\begin{equation}
 \{(M^{-1}\otimes I_n)+(I_n\otimes (-M))\}\bm{\tilde{x}}=\bm{c}'',\label{corollary1}
\end{equation}
where $\bm{c}'':=(M^{-1}\otimes I_n)\bm{c}'$.
Applying the inverse vec operator to (\ref{corollary1}), we have
\begin{align}
&{\rm vec}^{-1}\{\{(M^{-1}\otimes I_n)+(I_n\otimes (-M))\}\bm{\tilde{x}}\}={\rm vec}^{-1}\{\bm{c}''\}\nonumber\\
&\quad \Leftrightarrow -M\tilde{X}+\tilde{X}(M^{-1})^{\rm T}=Q'.\nonumber
\end{align}
This completes the proof. \hfill \ $\Box$
\par
\quad 

\section{Concluding remarks}
\quad In this paper, we showed that the T-congruence Sylvester equation can be transformed into the Lyapunov equation if the matrix $A$ or $B$ is nonsingular, and can be further transformed into the Sylvester equation if the matrix $A$ and $B$ are nonsingular. \\
\quad As applications, these results will lead to the following potential advantages: (1) simplification of the conditions for the unique solution of T-congruence Sylvester equation by using the condition for the Lyapunov equation, see, e.g., \cite{L_unique,S_unique}; (2)useful tools to find necessary and sufficient conditions for consistency of the T-congruence Sylvester equation; (3) efficient numerical solvers for the equation by using the numerical solvers of the Lyapunov equation (the Sylvester equation), see, e.g., \cite{S_unique,{Sylvester},Lyapunov}.\\
\quad In future work we will investigate whether there exists a relationship between the T-congruence Sylvester equation and the Lyapunov equation for the case where the matrix $A$ and $B$ are nonsingular or rectangular.
\begin{ack}
This work has been supported in part by JSPS KAKENHI Grant No. 26286088.
\end{ack}

\end{document}